\documentclass{ifacconf}

\usepackage{graphicx}
\usepackage{natbib}
\usepackage{amssymb,amsmath,dsfont}
\usepackage[english]{babel}
\usepackage{tikz,pgflibraryshapes,xcolor,xkeyval}
\usepackage{subfigure}
\usepackage{caption}

\newcommand {\etavec}{\boldsymbol{\eta}}

\newcommand {\varphivec}{{\boldsymbol{\varphi}}}

\newcommand {\col}{{\mathrm{col}}}

\newcommand {\alphavec}{\boldsymbol{\alpha}}

\newcommand {\xivec}{{\boldsymbol{\xi}}}
\newcommand {\lambdavec}{{\boldsymbol{\lambda}}}

\newcommand {\zetavec}{{\boldsymbol{\zeta}}}

\newcommand {\Psivec}{\mbox{\boldmath $\Psi$}}

\newcommand {\thetavec}{{\boldsymbol{\theta}}}

\newfont{\pseudocode}{cmtt10}

\newcommand{\bfx}{\mathbf{x}}
\newcommand{\bfl}{\boldsymbol{\ell}\,}

\newcommand{\bff}{\mathbf{f}}
\newcommand{\bfv}{\mathbf{v}}
\newcommand{\bfb}{\mathbf{b}}
\newcommand{\bfB}{\mathbf{B}}
\newcommand{\bfM}{\mathbf{M}}

\newcommand{\bfe}{\mathbf{e}}
\newcommand{\bfd}{\mathbf{d}}
\newcommand{\bfA}{\mathbf{A}}
\newcommand{\bfI}{\mathbf{I}}

\newcommand{\bfC}{\mathbf{C}}
\newcommand{\bfz}{\mathbf{z}}

\newcommand{\bfh}{\mathbf{h}}

\newcommand{\bfs}{\mathbf{s}}

\newcommand{\bfg}{\mathbf{g}}

\newcommand{\bfp}{\mathbf{p}}
\newcommand{\bfq}{\mathbf{q}}

\newcommand{\bfG}{\mathbf{G}}

\newcommand{\dist}{\mathrm{dist}}

\newcommand{\Real}{\mathbb{R}}

\newcommand{\Numbers}{\mathbb{Z}}

\newcommand{\norm}[1]{\left\Vert#1\right\Vert}

\renewcommand{\Real}{\mathbb{R}}

\renewcommand{\Real}{\mathds{R}}

\renewcommand{\Numbers}{\mathds{Z}}

\newcommand{\ra}{\rightarrow}

\newcommand{\Ecal}{\mathcal{E}}

\renewcommand{\col}[1]{\mathrm{col}\left(#1\right)}

\newcommand{\beq}{\begin{equation}}
\newcommand{\eeq}{\end{equation}}

\newcommand{\bbm}{\begin{bmatrix}}
\newcommand{\ebm}{\end{bmatrix}}

\newcommand{\bpm}{\begin{pmatrix}}
\newcommand{\epm}{\end{pmatrix}}

\newcommand{\bit}{\begin{itemize}}
\newcommand{\eit}{\end{itemize}}

\newcommand{\ben}{\begin{enumerate}}
\newcommand{\een}{\end{enumerate}}

\newcommand{\barr}{\begin{array}}
\newcommand{\earr}{\end{array}}

\renewcommand{\norm}[1]{\Vert#1\Vert}

\begin{document}

\begin{frontmatter}

\title{Adaptive observers for nonlinearly parameterized systems subjected to parametric constraints}

\author[First]{Ivan Yu. Tyukin}
\author[Second]{Pavel Rogachev}
\author[Third]{Henk Nijmeijer}

\address[First]{University of Leicester, Department of Mathematics, Leicester, United Kingdom (e-mail: I.Tyukin@le.ac.uk)}
\address[Second]{Saint-Petersburg State Electrotechnical University, Department of Automation and Control Processes, Saint-Petersburg, Russia (e-mail: rogapavel@yandex.ru)}
\address[Third]{Department of Mechanical Engineering, Eindhoven University of Technology, P.O. Box 513 5600 MB,  Eindhoven, The Netherlands, (h.nijmeijer@tue.nl)}

\begin{abstract}                % Abstract of not more than 250 words.

We consider the problem of adaptive observer design in the
settings when the system is allowed to be nonlinear in the
parameters, and furthermore they are to satisfy additional
feasibility constraints. A solution to the problem is proposed
that is based on the idea of universal observers and
non-uniform small-gain theorem. The procedure is illustrated
with an example.
\end{abstract}
\begin{keyword}
Observers, adaptive systems, state and parameter estimation, nonlinear parametrization.
\end{keyword}
\end{frontmatter}

\section{Introduction}

The problem of adaptive observer design has been in the focus of
considerable attention in the past few decades, see  e.g.
\cite{Bastin88},\cite{Marino90},\cite{MarinoTomei95},
\cite{MarinoTomei93}, \cite{Besancon:2000}. Available results
apply to a broad range of models, including to both linear and
non-linear in parameter systems of differential equations
\cite{Automatica:Farza:2009},\cite{SysConLett:Grip:2011}.

\begin{figure*}
\includegraphics[width=\textwidth]{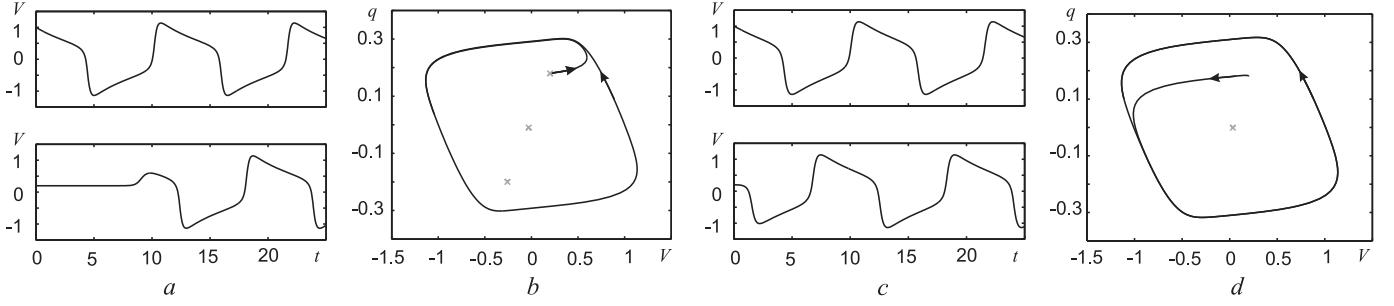}
\caption{Dynamics of (\ref{eq:Rowat_model}) for different parameter values.}\label{fig:motivation:1}
\end{figure*}

Despite this success there is still room for further
development. One particular direction is the case when the
model parameters are to satisfy additional constraints
corresponding to the feasibility regions in the parameter
space. Consider for instance the following system \cite{row}
\begin{equation}\label{eq:Rowat_model}
\begin{split}
\left(\begin{array}{l}
\dot{V}\\
\dot{q}
\end{array}\right)=&
\left(\begin{array}{cc}
       -\tau_m^{-1} & - \tau_m^{-1}\\
        \sigma_s\tau_s^{-1} & -{\tau}_s^{-1}
      \end{array}
\right)\left(\begin{array}{l} V \\ q \end{array}\right) \\
& \ \ \ \ \ \ \ \ \ \ \ \ \ \ \ \ + \left(\begin{array}{l} 1\\ 0\end{array}\right)\frac{A_f}{\tau_m} \tanh\left(\frac{\sigma_f}{A_f} V\right),
\end{split}
\end{equation}
where $\tau_m$, $\tau_s$, $\sigma_s$, $\sigma_f$, $A_f$ are
parameters. Suppose that variable $V$ is available for
measurement, and $q$ is not accessible for direct observation.
Equations (\ref{eq:Rowat_model}) is a modification of the
Van-der-Pol equations. Depending on the values of parameters,
model (\ref{eq:Rowat_model}) is known to exhibit a number of
oscillatory solutions with different qualitative dynamics. In
systems of this type there may exist a set of solutions for
which measured variables, e.g. $V$, look the same, and yet
their corresponding parameter values and qualitative underlying
dynamics are different. An example of this situation for (\ref{eq:Rowat_model}) is illustrated with  Fig. \ref{fig:motivation:1}. {\it Panel $a$} shows numerical approximation of trajectories  $V(t)$ in (\ref{eq:Rowat_model}) for $\tau_s=5$, $\tau_m=0.1666$, $A_f=1$, $\sigma_s=0.9$, and $\sigma_f=2$. Upper and lower plots correspond to solutions of (\ref{eq:Rowat_model}) passing through $q(0)=0$, $V(0)=1$ and $q(0)= 0.179756$, $V(0)=0.199729$, respectively. {\it Panel $b$} contains phase curves of the latter solutions. {\it Panel $c$} presents trajectories $V(t)$ corresponding to the same initial conditions as in panel $a$ but for different parameter values ($\sigma_s=1.1$, $\tau_s=6.1881$, $\tau_m=0.2062$, other parameters remain unchanged). Phase curve of the solution passing through $q(0)= 0.179756$, $V(0)=0.199729$ in for these modified parameter values is shown in {\it panel $d$}. Grey crosses in panels $b$ and $d$ indicate equilibria.
Note that trajectories $V(t)$ corresponding to solutions passing through $q(0)=0$, $V(0)=1$ (upper plots in Fig. \ref{fig:motivation:1}.$a$, Fig. \ref{fig:motivation:1}.$c$) look nearly identical despite that these trajectories are generated by systems (\ref{eq:Rowat_model})  with different parameter values and with different number of equilibria. The overall dynamics of these systems, as the other plots in  Fig. \ref{fig:motivation:1} suggest, are rather different.

The example above suggests that inferring true parameter values and estimating states from given output observations for systems capable of exhibiting multiple dynamic regimes is not always a well-posed problem. Lack of excitation and/or presence of unmodeled dynamics may result in the values of estimates that are not physically plausible. Hence additional care and precaution must be taken when dealing with such systems. Traditional methods of observer design do not always allow to account for parameter feasibility regions except, possibly, when the set of feasible parameters is convex \cite{SysConLett:Grip:2011,Kokotovich95}. Thus developing methods enabling to deal with nonlinearly parameterized systems subjected to rather general parametric constraints is needed.

Here we propose an adaptive observer design that satisfies these requirements. Our method is based on the idea of exploiting the advantage of combining
search-based optimization with that of direct optimization
routines. Observers of this sort were suggested and analyzed in
\cite{Tyukin:2013:Automatica}. In this work we generalize these
results by enabling the observer to account for parameter
feasibility regions.

The manuscript is organized as follows. In Section \ref{sec:problem_formulation} we describe the class of systems considered in the article and provide mathematical statement of the problem. Sections \ref{sec:observer_definition} and \ref{sec:main_result} present main results, Section \ref{sec:example} contains an illustrative example, and Section \ref{sec:conclusion} concludes the paper. Proofs of auxiliary technical statements are provided in the Appendix.

\section{Problem Formulation}\label{sec:problem_formulation}

We will deal with the following class of systems
\begin{equation}\label{eq:neural_model:3}
\begin{array}{l}
\dot{\bfx}=\bfA \bfx + \Psivec(t,\lambdavec,y)\thetavec+\bfg(t,\lambdavec,y,u)+\xivec(t),\\
y=\bfC^{T}\bfx, \ \bfA=\left(\begin{array}{cc}0 & \bfI_{n-1}\\ 0 &
0\end{array}\right), \ \bfC=\mathrm{col}(1,0,\dots,0),
\end{array}
\end{equation}
where
$\Psivec:\Real\times\Real^p\times\Real\rightarrow\Real^{n\times
m}$, $\Psivec,\bfg\in\mathcal{C}^1$, are bounded and Lipschitz
in $\lambdavec$, and $\dot{\Psivec}(t,\lambdavec,y(t))$,
$\dot{\bfg}(t,\lambdavec,y(t),u(t))$ are bounded;
$\lambdavec=\mathrm{col}(\lambda_1,\dots,\lambda_p)\in\Real^p$,
$\thetavec=\mathrm{col}(\theta_1,\dots,\theta_m)\in\Real^{m}$
are {\it unknown parameters}, and $u:\Real\rightarrow\Real$,
$u\in\mathcal{C}^1$ is the input. We assume that the values of
$\lambdavec$, $\thetavec$ belong to the hypercubes
$\Omega_\lambda\subset\Real^p$, $\Omega_\theta\subset\Real^{m}$
with known bounds:
$\theta_{i}\in[\theta_{i,\min},\theta_{i,\max}]$,
$\lambda_j\in[\lambda_{j,\min},\lambda_{j,\max}]$. The function
$\xivec\in\mathcal{C}^0:\Real\rightarrow\Real^n$ represents
{\it unmodeled dynamics} and is supposed to be unknown yet
bounded: \beq\label{eq:xi_bound}
    \|\xivec(t)\|\leq \Delta_\xi, \ \Delta_\xi\in\Real_{\geq 0}, \ \forall \ t.
\eeq

Since not all systems (\ref{eq:neural_model:3})
may be identifiable, we will require parameter reconstruction only
up to a set $\mathcal{E}(\lambdavec,\thetavec)$ of indistinguishable parameters that will be specified later. In addition, we suppose
that, for the observed trajectories of
(\ref{eq:neural_model:3}), true values of parameters $\lambdavec,\thetavec$
satisfy the following constraint:
\[
{\pi}(\thetavec,\lambdavec)=0,
\]
where $\pi:\Real^m\times\Real^p\rightarrow\Real$ is a Lipschitz
function. The problem is to determine an auxiliary system, i.e.
an {\it adaptive observer}: $\dot{\bfq}=\bff(t,y,u,\bfq)$,
$\bff:\Real\times\Real\times\Real\times\Real^{q}\rightarrow\Real^q$,
$\bfq(t_0)=\bfq_0$, $\bfq\in\Real^q$, and functions
$\bfh_{\theta}:\Real^q\rightarrow\Real^m$,
$\bfh_{\lambda}:\Real^q\rightarrow \Real^p$,
$\bfh_x:\Real^q\rightarrow\Real^n$ such that for some given
$\bfq_0$,  known functions $r_1,r_2,r_3\in\mathcal{K}$, and all
$t_0\in\Real$, the following requirements hold  for the
observer:
\begin{eqnarray}\label{eq:goal}
     & &\limsup_{t\rightarrow\infty}  \|\bfh_{x}(\bfq(t,\bfq_0))-\bfx(t,\lambdavec,\thetavec,\bfx_0)\|\leq r_1(\Delta_{\xi}) \\
    %&\limsup_{t\rightarrow\infty} \|\bfh_{\theta}(\bfq(t,\bfq_0))-\thetavec\|\leq r_2(\Delta_{\xi})\\
    & &\limsup_{t\rightarrow\infty}   \
        \dist\left(\left(\begin{array}{c}\bfh_{\lambda}(\bfq(t,\bfq_0))\\
                                \bfh_{\theta}(\bfq(t,\bfq_0))
    \end{array}\right),\mathcal{E}(\lambdavec,\thetavec)\right)\leq r_2(\Delta_{\xi}),\nonumber\\
    & & \|\pi(\bfh_{\theta}(\bfq(t,\bfq_0)),\bfh_{\lambda}(\bfq(t,\bfq_0)))\|\leq r_3(\Delta_\xi),
    \nonumber
\end{eqnarray}
where $\Delta_\xi$ is defined in (\ref{eq:xi_bound}).

\section{Observer Definition}\label{sec:observer_definition}

As an observer candidate for (\ref{eq:neural_model:3}) we
propose the following:
\begin{eqnarray}\label{eq:linear_par_observer:modified}
&&\dot{\bfM}=(\bfA-\bfB\bfC^{T}\bfA)\bfM+(\bfI_n-\bfB\bfC^{T})\Psivec(t,\hat{\lambdavec},y),\nonumber\\
&&\dot{\hat{\zetavec}}=\bfA\hat{\zetavec}+\bfl(\bfC^{T}\hat{\zetavec}-y)+\bfB{\varphivec}^{T}(t,\hat{\lambdavec},y,[\hat\lambdavec,y])\hat{\thetavec}\nonumber\\
&& \ \ \ \ \ \ \ \ \ \ \ \ \ \ + \bfg(t,\hat{\lambdavec},y,u),  \\
&&\dot{\hat{\thetavec}}=-\gamma_\theta (\bfC^{T}\hat{\zetavec}-y){\varphivec}(t,\hat{\lambdavec},y,[\hat\lambdavec,y]), \ \gamma_\theta\in\Real_{>0},\nonumber\\
&& \hat{\bfx}=\hat{\zetavec}+\bfM\hat{\thetavec}\nonumber,
\end{eqnarray}
where $\bfB=\col{1,b_1,\dots,b_{n-1}}$ is a vector such that
the polynomial $s^{n-1}+b_1 s^{n-2}+\dots+b_{n-1}$ is Hurwitz,
$\bfM(t_0)=0$,
\[
{\varphivec}^{T}(t,\hat{\lambdavec},y,[\hat\lambdavec,y])=\bfC^{T}\bfA \bfM(t,[\hat{\lambdavec},y]) +\bfC^{T}\Psivec(t,\hat{\lambdavec},y),
\]
and $\hat\lambdavec$ is defined as
\begin{equation}\label{eq:nonlinear_par_observer}
\begin{split}
                &\left\{\begin{array}{ll}
                    \dot{{s}}_{2j-1}&=\gamma \sigma(\|y-\hat{y}\|_\varepsilon+\|\pi(\hat{\thetavec},\hat{\lambdavec})\|_\varepsilon) \cdot\omega_j  \cdot ({s}_{2j-1} \\
                    & - {s}_{2j}- {s}_{2j-1}({s}_{2j-1}^2+{s}_{2j}^2))
                    \\
                    \dot{{s}}_{2j}&=\gamma \sigma(\|y-\hat{y}\|_\varepsilon+\|\pi(\hat{\thetavec},\hat{\lambdavec})\|_\varepsilon)\cdot\omega_j \cdot ({s}_{2j-1} \\
                    &+ {s}_{2j}- {s}_{2j}({s}_{2j-1}^2+{s}_{2j}^2) )
                    \\
                    \hat{\lambda}_j&=\beta_j(\bfs), \ j=\{1,\dots,p\},
                 \end{array}\right.\\
                 \beta_j(\bfs)&=\lambda_{j,\min}
                    +\frac{\lambda_{j,\max}-\lambda_{j,\min}}{2}({s}_{2j-1}+1)
            \end{split}
        \end{equation}
\begin{equation}\label{eq:initial_conditions}
     \bfs_0=\bfs(t_0): \ {s}_{2j-1}^2(t_0)+{s}_{2j}^2(t_0)=1,
\end{equation}
where $\gamma\in\Real_{>0}$,  $\sigma:\Real_{\geq
0}\rightarrow\Real_{\geq 0}$ is a bounded Lipschitz function:
\begin{equation}%\label{eq:generator_poisson:3}
\begin{split}
\exists  \ D_\sigma, \ M_{\sigma}\in\Real_{>0}: &\\
  \sigma(\upsilon)\leq M_{\sigma},&  \ \sigma(\upsilon)\leq D_\sigma \upsilon \  \ \forall \  \upsilon\geq 0
\end{split}\nonumber
\end{equation}
such that $\sigma(\upsilon)>0$ for $\upsilon>0$, and
$\sigma(0)=0$.  Parameters $\omega_j\in\Real_{>0}$ in
(\ref{eq:nonlinear_par_observer}) are supposed to be {\it
rationally-independent}:
\begin{equation}\label{eq:rational_independence}
    \sum_{j=1}^p \omega_j k_j\neq 0, \ \forall \ k_j\in \Numbers.
\end{equation}

\section{Main Result}\label{sec:main_result}

As it is often the case in the domain of adaptive observer
design, state and parameter reconstruction of the observer is
subjected to some form of persistency of excitation (PE). Insufficient excitation have detrimental effect on the quality and robustness of estimation (see examples in Fig. \ref{fig:motivation:1}). Here
we employ the following modifications of the PE conditions:
\begin{defn}[\cite{Lorea_2002}]\label{defn:uniform_pe}
       A function $\alphavec:\Real_{\geq t_0} \times \Omega_\lambda \ra \Real^p$ is said to be $\lambda$-Uniformly Persistently Exciting ($\lambda$-UPE with $T$, $\mu$), denoted by $\alphavec(t,\lambdavec)\in\lambda{\rm UPE}(T,\mu)$, if there exist $T,\mu\in\Real_{>0}$:
    \beq\label{eq:PE_linear_uniform}
        \int_{t}^{t+T} \alphavec(\tau,\lambdavec)\alphavec^T(\tau,\lambdavec){\rm d}\tau \geq \mu
        \bfI, \ \forall \ t\geq t_0, \ \lambdavec\in
        \Omega_\lambda.
    \eeq
\end{defn}
Uniform persistency of excitation requires existence of
$\mu\in\Real_{>0}$ in (\ref{eq:PE_linear_uniform}) that is
independent on $\lambdavec$ for all
$\lambdavec\in\Omega_\lambda$.
\begin{defn}\label{defn:nonlinear_pe} Let $\mathcal{E}$ be
a set-valued map defined on $\mathcal{D}\subset\Real^d$  and
associating a subset of $\mathcal{D}$ to every
$\bfp\in\mathcal{D}$.
   A function $\alphavec:\Real \times \mathcal{D}\times\mathcal{D} \ra \Real^k$ is
    said to be weakly Nonlinearly Persistently Exciting in $\bfp$ wrt $\mathcal{E}$ (wNPE with $L,\beta,\mathcal{E}$ ),
    denoted by $\alphavec(t,\bfp,\bfp')\in{\rm wNPE}(L,\beta,\mathcal{E})$, if there exist $L\in\Real_{>0}$, $t_1\geq t_0$, and $\beta\in\mathcal{K}_{\infty}$:
    \beq\label{eq:nonlinear_pe}
    \begin{split}
     & \forall \ t\geq t_1, \ \bfp,\bfp'\in\mathcal{D}  \ \exists \  t'\in[t,t+L]:  \\
     & \ \ \ \ \ \ \ \ \norm{\alphavec(t',\bfp,\bfp')}
     \geq \beta \left({\rm dist}(\Ecal(\bfp),\bfp^\prime)\right).
    \end{split}
    \eeq
\end{defn}
Finally, consider
\[
\begin{split}
&\etavec(t,\lambdavec,\thetavec,\lambdavec',\thetavec')=\varphivec^{T}(t,\lambdavec',y(t),[\lambdavec',y])(\thetavec'-\thetavec)+\\
&\bfC^{T}(\Psivec(t,\lambdavec',y(t))-\Psivec(t,\lambdavec,y(t))\thetavec + g_1(t,\lambdavec',y(t),u(t))-\\
& g_1(t,\lambdavec,y(t),u(t))+ q(t,\lambdavec',\lambdavec,\thetavec),
\end{split}
\]
 $q(t,\lambdavec',\lambdavec,\thetavec)=\tilde{\bfC}\bfz$,
$\dot{\bfz}=\Lambda \bfz + \bfG(
(\Psi(t,\lambdavec',y(t))-\Psi(t,\lambdavec,y(t)))\thetavec+\bfg(t,\lambdavec',y(t),u(t))-\bfg(t,\lambdavec,y(t),u(t)))$,
$\bfz(t_0)=0$,  and $\tilde{\bfC}$, $\Lambda$, $\bfG$ are
defined as
\[
\begin{split}
&\tilde{\bfC}=\mathrm{col}(1,0,\dots,0), \tilde{\bfC} \in \Real^{n-1}, \
\Lambda=\left(\begin{array}{ccc}
-\bfb \  & \begin{array}{c}\vdots\\ \vdots \end{array} &
\begin{array}{c} I_{n-2}\\ 0
\end{array}\end{array}\right), \\
&\bfG=\left(\begin{array}{cc} - {\bfb} & I_{n-1} \end{array}\right), \
\bfb=(b_1,\dots,b_{n-1})^{T}.
\end{split}
\]
The following can now be stated:
\begin{thm}\label{cor:general} Consider \eqref{eq:neural_model:3},  \eqref{eq:linear_par_observer:modified}, \eqref{eq:nonlinear_par_observer}--\eqref{eq:rational_independence}. Suppose that the restriction of the function $\varphivec(t,\lambdavec,y(t),[\lambdavec,y])$ in (\ref{eq:linear_par_observer:modified}) on
    $\Real\times\Omega_{\lambda}$ is $\lambda$-uniformly
    persistently exiting, and the function
$\alphavec(t,(\lambdavec,\thetavec),(\lambdavec',\thetavec'))=\etavec(t,\lambdavec,\thetavec,\lambdavec',\thetavec')$
is weakly nonlinearly persistently exciting in
$(\lambdavec,\thetavec)$ wrt to the map $\mathcal{E}$:
\[
\begin{split}
&\mathcal{E}(\lambdavec,\thetavec)=\{(\lambdavec',\thetavec'), \ \lambdavec'\in\Real^p, \ \thetavec'\in\Real^m| \bfB (\thetavec'-\thetavec)^{T}\cdot \\
& \varphivec(t,\lambdavec',y(t),[\lambdavec',y])+ (\Psivec(t,\lambdavec',y(t))-\Psivec(t,\lambdavec,y(t)))\thetavec\\
&+\bfg(t,\lambdavec',y(t),u(t))-\bfg(t,\lambdavec,y(t),u(t))=0,\ \forall \ t\geq t_0\}.
\end{split}
\]
Then there exist a constant $\bar \gamma \in \Real_{>0}$ and
functions $r_1,r_2,r_3,r_4\in\mathcal{K}$ such that if $\gamma,
\varepsilon$ are the corresponding parameters of
(\ref{eq:nonlinear_par_observer}), and $\gamma \in (0, \bar
\gamma)$, $\varepsilon > r_1(\Delta_{\xi})$,  then
\begin{align}
        &\limsup_{t\ra\infty}
        \dist\left(\left(\begin{array}{c}\hat\lambdavec(t)\\
                                \hat\thetavec(t)
    \end{array}\right),\mathcal{E}(\lambdavec,\thetavec)\right)\leq r_2(\varepsilon). \label{eq:thm2}
\end{align}
\begin{align}
        &\limsup_{t\ra\infty}\norm{\hat \bfx(t)-\bfx(t)} \leq r_3(\varepsilon), \label{eq:thm1}
\end{align}
\begin{align}
        &\limsup_{t\ra\infty}
        \|\pi(\hat{\thetavec}(t),\hat{\lambdavec}(t))\| \leq r_4(\varepsilon). \label{eq:thm3}
\end{align}
\end{thm}
{\it Proof of Theorem \ref{cor:general}}. The idea behind the
proof of the theorem is similar to that of presented in
\cite{Tyukin:2013:Automatica} and, taking space limitations into
account, we present just a sketch here. Note that the first row
of $\bfM$ in the observer system is zero for all $t\geq t_0$, and
that $\hat{y}=\bfC^{T}\hat{\bfx}=\bfC^{T}\hat{\zetavec}$.
Moreover, since $\Psivec(t,\hat{\lambdavec},y)$ is bounded and
Lipschitz in $\hat{\lambdavec}$, the variables
$\bfM(t,[\hat{\lambdavec},y]), \dot{\bfM}(t,[\hat{\lambdavec},y])$
are also bounded, and ${\bfM}(t,[\hat\lambdavec,y])$ is Lipschitz
in $\hat{\lambdavec}$ for $\hat\lambdavec=\mathrm{const}$. Let
$\zetavec=\bfx-\bfM\thetavec$, then
\[
\begin{split}
{\dot\zetavec}=&\bfA \zetavec
+ \bfB\varphivec(t,\hat\lambdavec,y,[\hat{\lambdavec},y])\thetavec +
(\Psivec(t,\lambdavec,y,u)\\
&-\Psivec(t,\hat\lambdavec,y,u))\thetavec
+ \bfg(t,\lambdavec,y,u)+\xivec(t).
\end{split}
\]
Dynamics of (\ref{eq:neural_model:3}),
(\ref{eq:linear_par_observer:modified}) in the coordinates
$\bfe_1=\hat{\zetavec}-\bfx+\bfM\thetavec$,
$\bfe_2=\hat{\thetavec}-\thetavec$ is
\begin{eqnarray}\label{eq:error_sys:modified}
&& \bpm \dot \bfe_1 \\ \dot \bfe_2 \epm
        =
        \bpm
            \bfA+\bfl \bfC^T   &   \bfB \varphivec^T(t,\hat \lambdavec,y,[\hat \lambdavec,y]) \\
            -\gamma_{\theta} \varphivec(t,\hat \lambdavec,y,[\hat \lambdavec,y])\bfC^T & 0
        \epm \nonumber\\
&& \ \ \ \ \ \ \ \ \ \ \ \ \ \ \ \times\bpm
            \bfe_1 \\ \bfe_2
        \epm
         +
        \bpm
            \tilde{\bfv}(t,\hat \lambdavec,\lambdavec, y,u) \\ 0
        \epm
\end{eqnarray}
where
\[
\begin{split}
\tilde{\bfv}(t,\hat \lambdavec,\lambdavec,
y,u)&=(\Psivec(t,\hat
\lambdavec,y)-\Psivec(t,{\lambdavec},y))\thetavec\\
&+ \bfg(t,\hat
\lambdavec,y,u)-\bfg(t,\lambdavec,y,u)-\xivec(t).
\end{split}
\]
Since the
pair $\bfA$, $\bfC$ is observable one can always find an $\bfl$
so that dynamics of the homogeneous part of
(\ref{eq:error_sys:modified}) is exponentially stable subject
to persistency of excitation of $\varphivec(t,\hat
\lambdavec,y,[\hat \lambdavec,y])$. The latter condition is
ensured by picking the values of $\gamma$ sufficiently small
(see \cite{Tyukin:2013:Automatica}).

Recall that  the function $\pi$ is Lipschitz and that
$\pi(\thetavec,\lambdavec)=0$. Therefore
\begin{equation}\label{eq:constraint:aux}
\begin{split}
\|\pi(\hat{\thetavec},\hat{\lambdavec})\|&=\|\pi(\hat\thetavec,\hat\lambdavec)-\pi(\thetavec,\lambdavec)\|
\leq L_1 \|\hat\thetavec-{\hat\thetavec}\|\\
&+L_2 \|\hat\lambdavec-\lambdavec\|.
\end{split}
\end{equation}
Furthermore, notice that
\begin{enumerate}
\item $|a|\leq |b| \Rightarrow \|a\|_\varepsilon\leq \|b\|_\varepsilon$
\item $\||a|+|b|\|_\varepsilon\leq \|a\|_\varepsilon+|b| \ \forall \ a,b\in\Real$
\item $ \|L a\|_\varepsilon = |L| \|a\|_{\frac{\varepsilon}{|L|}},  \ a \in\Real$
\item $\|a\|_\varepsilon + \|b\|_\varepsilon \leq 2 \|a,b\|_\varepsilon$, where $\|a,b\|_\varepsilon = \sqrt{\|a\|_\varepsilon^2+\|b\|_\varepsilon^2}$.
\end{enumerate}
Hence, using (\ref{eq:constraint:aux}) and applying properties
1--4 above to
$\|\pi(\hat\thetavec,\hat\lambdavec)\|_\varepsilon$ we obtain:
\[
\begin{split}
\|\pi(\hat{\thetavec},\hat{\lambdavec})\|_\varepsilon&\leq \|L_1\|\hat\thetavec-\thetavec\| \|_\varepsilon + L_2 \|\hat\lambdavec-\lambdavec\| = L_1 \|\hat\thetavec-\thetavec\|_\frac{\varepsilon}{L_1}\\
&+ L_2 \|\hat\lambdavec-\lambdavec\|.
\end{split}
\]
Therefore, since $\|\bfC\|=1$, $\|\bfC^{T}(\bfx-\hat{\bfx})\|_\varepsilon + \|\pi(\hat\thetavec,\hat\lambdavec)\|_\varepsilon \leq \|\hat\bfx-\bfx\|_\varepsilon + L_1 \|\hat\thetavec-\thetavec\|_\frac{\varepsilon}{L_1}+ L_2 \|\hat\lambdavec-\lambdavec\|$.
Denoting $\epsilon=\frac{\varepsilon}{L_1}$ results
in
\begin{equation}\label{eq:constraint:aux:estimate}
\begin{split}
& \|\bfC^{T}(\bfx-\hat{\bfx})\|_\varepsilon + \|\pi(\hat\thetavec,\hat\lambdavec)\|_\varepsilon\leq \|\hat\bfx-\bfx\|_\epsilon\\
&+L_1\|\hat\thetavec-\thetavec\|_\epsilon+L_2\|\hat\lambdavec-\lambdavec\|.
\end{split}
\end{equation}
Finally, taking (\ref{eq:constraint:aux:estimate}) into
account, and repeating the argument provided in
\cite{Tyukin:2013:Automatica} in which Lemma 8 is replaced with
Lemma \ref{lem:bound} below (see Appendix for the proof) one
can complete the proof of the theorem.
\begin{lem}\label{lem:bound}
Consider a system governed by the following set of equations
\begin{eqnarray} \label{eq:interconnection}
        & &\norm{\bfx(t)} \leq \beta(t-t_0)\norm{\bfx(t_0)}+c \norm{h(\tau)}_{\infty,[t_0,t]} + \Delta,\\
       & &- \int_{t_0}^t \gamma_0(\norm{\bfx(\tau)+\bfd(\tau)}_\varepsilon+M|h(\tau)|)d\tau \nonumber \\
       & & \ \ \ \ \ \ \ \ \ \ \ \ \ \  \ \ \ \ \ \ \ \ \ \ \ \ \ \ \ \ \leq h(t) - h(t_0) \leq 0, \  \forall \ t\geq t_0, \nonumber
 \end{eqnarray}
where $\bfx:\Real\rightarrow\Real^n$, $h:\Real\rightarrow\Real$
are trajectories reflecting the evolution of the system's
state, $\bfd:\Real\rightarrow\Real^n$,
$\|\bfd(\tau)\|_{\infty,[t_0,\infty)}\leq \Delta_d$ is a
continuous and bounded function on $[t_0,\infty)$,
$\beta(\cdot)$ is a strictly monotonically decreasing function
with, $\beta(0)\geq 1$, $\lim_{s\ra\infty}\beta(s)=0$;
$c,\Delta\in \Real_{>0}$, and
$\gamma_0:\Real\rightarrow\Real_{\geq 0}$:
\begin{equation}\label{eq:gamma_0}
|\gamma_0(s)|\leq D_\gamma |s|.
\end{equation}
Then trajectories $\bfx(t)$, $h(t)$ in
(\ref{eq:interconnection}) are bounded in forward time, for
$t\geq t_0$, provided that the following conditions hold for
some $d\in(0,1)$, $\kappa\in(1,\infty)$:
\begin{eqnarray}
        \varepsilon&\geq& \Delta \left(1+\beta(0)\frac{\kappa}{\kappa - d}\right)+\Delta_d, \label{eq:epsilon_choice}\\
       D_\gamma &\leq&\frac{\kappa-1}{\kappa} \left[\beta^{-1}\left(\frac{d}{\kappa}\right)\right]^{-1} \times \label{eq:gamma_choice}\\
&  &  \times  \frac{h(t_0) }{ \beta(0)\|\bfx(t_0)\|+ |h(t_0)|(c(1+ \kappa \beta(0)/(1-d))+M)}.\nonumber
\end{eqnarray}
\end{lem}
$\square$

\section{Example}\label{sec:example}

Consider system (\ref{eq:Rowat_model}). Let us first transform (\ref{eq:Rowat_model}) into (\ref{eq:neural_model:3}). Applying the following coordinate transformation $
\left(
\begin{array}{c}
  V' \\
  q'
\end{array}
\right)
=\left(
\begin{array}{cc}
    1 & 0 \\
    0 & -\tau_m^{-1} \\
\end{array}
\right)
\left(
\begin{array}{c}
  V \\
  q
\end{array}
\right)$ we get:
$$
\begin{aligned}
\frac {dV'} {dt} &= - \tau_m^{-1} {V'} + \tau_m^{-1} {A_f \tanh [(\sigma_f/A_f)V']} + q',\\
\frac {dq'} {dt} &= - {\tau_s}^{-1} {q'} - {(\tau_s \tau_m)}^{-1} {\sigma_s} V',
\end{aligned}
$$
and using an additional change of coordinates $
\left(
\begin{array}{c}
  x_1 \\
  x_2
\end{array}
\right)
=\left(
\begin{array}{cc}
    1 & 0 \\
    \tau_s^{-1} & 1 \\
\end{array}
\right)
\left(
\begin{array}{c}
  V' \\
  q'
\end{array}
\right)$
we arrive at:
$$
\begin{aligned}
\dot{x}_1 =& x_2 - \left({\tau_m^{-1}} + {\tau_s^{-1}} \right) x_1 + \tau_m^{-1} A_f \tanh [(\sigma_f/A_f)x_1],\\
\dot{x}_2 =& -(\tau_s\tau_m)^{-1}\left(1 + {\sigma_s}\right) x_1 \\
&+ {(\tau_s\tau_m)^{-1}} A_f \tanh [(\sigma_f/A_f)x_1].
\end{aligned}
$$
Noticing that $y=C^T \mathbf{x}=x_1$ and denoting
$\theta_1=-\frac1{\tau_m} - \frac1{\tau_s},\quad \theta_2=\frac1{\tau_m} A_f,\quad \theta_3=\frac{-1-\sigma_s}{\tau_s\tau_m}, \quad \theta_4= \frac1{\tau_s\tau_m} A_f,\quad \lambda=\frac{\sigma_f}{A_f}$, we obtain
$$
\begin{aligned}
\dot{x_1} &= x_2 + \theta_1 y + \theta_2 \tanh [\lambda y],\\
\dot{x_2} &= \theta_3 y + \theta_4 \tanh [\lambda y],\\
y&=x_1.
\end{aligned}
$$
The above equations in the form (\ref{eq:neural_model:3})
where
$$
\mathbf{A}=\left(
  \begin{array}{cc}
    0 & 1 \\
    0 & 0 \\
  \end{array}
\right),\quad
\mathbf{\Psi}=\left(
  \begin{array}{cccc}
    y & \tanh [\lambda y] & 0 & 0 \\
    0 & 0 & y & \tanh [\lambda y] \\
  \end{array}
\right),
$$
$$
\thetavec=(\theta_1, \theta_2, \theta_3,  \theta_4), \
\mathbf{g}=\left(
  \begin{array}{c}
    0 \\
    0 \\
  \end{array}
\right)
,\quad
\mathbf{C}^T=\left(
     \begin{array}{cc}
       1 & 0 \\
     \end{array}
   \right).
$$
Parameters of the observer (\ref{eq:linear_par_observer:modified}), (\ref{eq:nonlinear_par_observer}) were set as follows:
$$
\mathbf{B}=\left(
  \begin{array}{c}
    1 \\
    1 \\
  \end{array}
\right),\quad
\mathbf{l}=\left(
  \begin{array}{c}
    -2 \\
    -1 \\
  \end{array}
\right),\quad
\gamma_{\mathbf{\theta}} = 4, \ \gamma=0.002.
$$

It is clear that original state parameter values of (\ref{eq:Rowat_model}) can be recovered from $\bfx$, $\thetavec$, and $\lambda$ in accordance with the inverse transform:
$$
\tau_s={\theta_2}/{\theta_4},\quad \tau_m=-1/({\theta_1+{1}/{\tau_s}}),\quad A_f=\theta_2 \tau_m,
$$
$$
\sigma_s=-\theta_3 \tau_s \tau_m -1, \ \sigma_f=\lambda A_f, \ \mathbf{x}=\mathbf{\zetavec}+\bfM\mathbf{\thetavec};
$$
$$
\left(
  \begin{array}{c}
    V \\
    q \\
  \end{array}
\right)=
\left(
  \begin{array}{cc}
    1 & 0 \\
    {\tau_s}{\tau_m^{-1}} & -{\tau_m^{-1}} \\
  \end{array}
\right)
\left(
  \begin{array}{c}
    x_1 \\
    x_2 \\
  \end{array}
\right).
$$

True parameter value and initial conditions were set to:
$$
\tau_s=5,\quad \tau_m=0.1666,\quad A_f=1, \quad \sigma_s=0.8, \quad \sigma_f=2;
$$
$$
\mathbf{\theta}=(-6.2, 6, -2.16, 1.2), \
\lambda=2,  \
\left(
  \begin{array}{c}
    V_0 \\
    q_0 \\
  \end{array}
\right)=
\left(
  \begin{array}{c}
    1 \\
    0 \\
  \end{array}
\right).
$$

\begin{figure}
\centering
\includegraphics[width=0.85\linewidth]{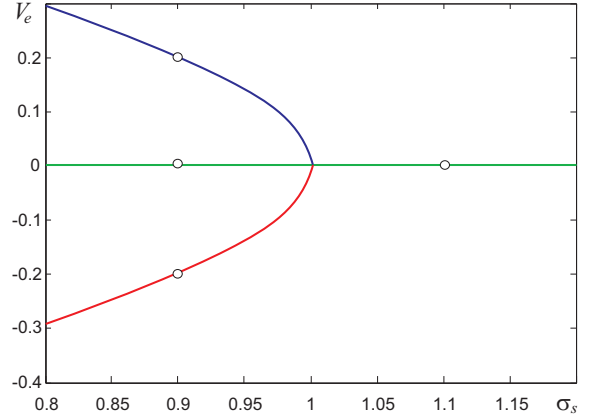}
\caption{Bifurcation diagram of system (\ref{eq:Rowat_model}).  }\label{fig:bifurcation_diagram}
\end{figure}

Depending on specific parameter values, the model may have $1$ or $3$ distinct equilibria in the system state space. One of these equilibria is always at the origin, $(0;0)$, and the other two are located in the first and third quadrant, respectively. The corresponding bifurcation diagram derived at $A_f=1$, $\sigma_f=2$ and $\sigma_s\in[0.8,1.2]$ is shown in Fig. \ref{fig:bifurcation_diagram}. Blue, green and red curves show the values of the $V$-coordinate, $V_e$, of the equilibria for various values of $\sigma_s$.  White circles correspond to solutions shown in panels $b$, $d$ in Fig. \ref{fig:motivation:1}.
 If $\sigma_s\leq\sigma_f-1$ the non-zero equilibria disappear. In our particular case the system has the following three equilibria: $(0;0)$, $(0.292;0.2336)$, $(-0.292;-0.2336)$. The first equilibrium is a saddle point (one positive and one negative real eigenvalues), and the other two are unstable nodes (both eigenvalues are real and positive).

We supposed that information about qualitative dynamics of the original system has been made available in the form of the number of equilibria and their qualitative characterizations (saddle point, stable/unstable nodes). In particular, we assumed that the following is known:
\begin{itemize}
  \item The total number of distinct equilibria in the system is $3$;
  \item Eigenvalues corresponding to zero equilibrium position are real (discriminant $> 0$);
  \item Signs of the eigenvalues at zero equilibrium are different.
\end{itemize}

The first restriction is satisfied if
\begin{equation}\label{eq:constraint:1}
\sigma_s>\sigma_f-1.
\end{equation}
Consider the other two remaining conditions.
Characteristic polynomial of the Jacobian of the right-hand side of (\ref{eq:Rowat_model}) is:
$s^2+\left((1-\sigma_f){\tau_m^{-1}}+{\tau_s^{-1}}\right)s+({1-\sigma_f+\sigma_s})(\tau_s\tau_m)^{-1}$.
The second condition is therefore equivalent to that the discriminant $D$ of this quadratic is positive:
\begin{equation}\label{eq:constraint:2}
D = \left(\frac{1-\sigma_f}{\tau_m}+\frac1{\tau_s}\right)^2 - 4 \left(\frac{1-\sigma_f+\sigma_s}{\tau_s\tau_m}\right)>0.
\end{equation}
The third and the final condition can now be expressed as:
\begin{equation}\label{eq:constraint:3}
\left|\frac{1-\sigma_f}{\tau_m}+\frac1{\tau_s}\right|<\sqrt D.
\end{equation}
Parametric constraints (\ref{eq:constraint:1})--(\ref{eq:constraint:3}) are in the form of strict inequalities. The proposed observer design, however, requires that parametric constraints are formulated in the form of equalities. For this purpose we introduce:
\[
\begin{split}
f_1(\thetavec,\lambda)=&\left|(1+\tanh(-u_1(\thetavec,\lambda)-3))/2\right|_\epsilon, \\  u_1(\thetavec,\lambda) =& \sigma_s(\thetavec,\lambda)-\sigma_f(\thetavec,\lambda)+1,\\
f_2(\thetavec,\lambda)=&\left|(1+\tanh(-u_2(\thetavec,\lambda)-3))/2\right|_\epsilon, \\  u_2(\thetavec,\lambda) =& D(\thetavec,\lambda),\\
f_3(\thetavec,\lambda)=&\left|(1+\tanh(-u_3(\thetavec,\lambda)-3))/2\right|_\epsilon, \\  u_3(\thetavec,\lambda)  =& \sqrt{D(\thetavec,\lambda)} - \left|\frac{1-\sigma_f(\thetavec,\lambda)}{\tau_m(\thetavec,\lambda)}+\frac1{\tau_s(\thetavec,\lambda)}\right|,
\end{split}
\]
where $\epsilon<(1+\tanh(-3))/2$, and define ${\pi}(\thetavec,\lambda)$ as follows:
\[
{\pi}(\thetavec,\lambda) = \sum_{i=1}^3 f_i(\thetavec,\lambda).
\]
Note that points $\thetavec,\lambda$ for which $\pi(\thetavec,\lambda)=0$ automatically satisfy conditions (\ref{eq:constraint:1})--(\ref{eq:constraint:3}).

\begin{figure*}[t]
\includegraphics[width=1\linewidth]{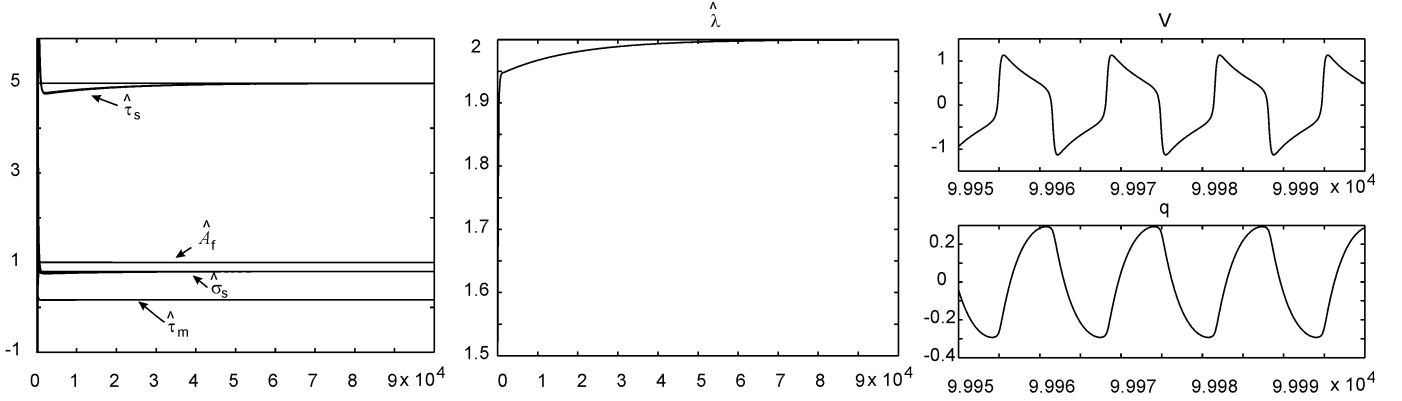}
\caption{Left panel: estimates of parameters  $\tau_m, \tau_s, \sigma_s, A_f$ as functions of $t$. Middle panel: estimate of  ${\lambda}$. Right panel: true values of $V$ and $q$ and their estimates in the end of the simulation.}\label{fig:results:ex1}
\end{figure*}

Results of numerical simulation of (\ref{eq:Rowat_model}) coupled with the observer are shown in Fig. \ref{fig:results:ex1}. As one can see from this figure, state and parameter estimates asymptotically approach their true values. In addition, we observed that including parametric constraints terms $\pi(\cdot,\cdot)$ into the observer equations not only helped to remove potential ill-conditioning of the problem but also it improved overall convergence times as compared to observers in which terms $\pi(\cdot,\cdot)$ have been dropped.

\section{Conclusion}\label{sec:conclusion}

In this work we provided an extension of the observer design for systems nonlinear in parameter \cite{Tyukin:2013:Automatica}. The extension allows to account for parametric constraints in the form of equalities. The approach presented enables inclusion of additional information about location of true parameter values and as such helps to improve robustness of the estimation procedure.

%%%%%%%%%%%%%%%%%%%%%%%%%%%%%%%%%%%%

%\bibliographystyle{ifacconf-harvard}
\bibliography{adaptive_observer_noncanonic_SCL}

\section*{Appendix}

{\it Proof of Lemma \ref{lem:bound}.} According to condition
(\ref{eq:gamma_choice}) we can conclude that that $h(t_0)\geq
0$. Let us introduce a strictly decreasing sequence: $
\{\sigma_i\},  \ i=0,1,\dots,  \ \sigma_i=(1/\kappa)^i, \
\kappa\in(1,\infty)$. Further, let $ \{t_i\}, \ i=1,\dots, \
t_1<t_2<\cdots<t_n<\cdots$ be an ordered infinite sequence of
time instants such that
\begin{equation}\label{eq:h_i}
 h(t_i)=\sigma_i h(t_0).
\end{equation}
If the latter assumption does not hold then one can immediately
conclude that $h(t)$ is bounded from below by $0$ for $t\geq
t_0$; it is  also bounded from above by $h(t_0)$, for $t\geq
t_0$. Hence $h(t)$ is bounded for all $t\geq t_0$. Moreover, in
accordance with (\ref{eq:interconnection}), trajectory
$\bfx(t)$ is bounded for all $t\geq t_0$, and nothing remains
to be proven.

We wish to show that if (\ref{eq:epsilon_choice}),
(\ref{eq:gamma_choice}), and (\ref{eq:h_i}) then
\begin{equation}\label{eq:limit_h}
h(t)\rightarrow 0 \Rightarrow t\rightarrow\infty.
\end{equation}
In order to do so consider the time differences
$T_i=t_i-t_{i-1}$. It is clear from (\ref{eq:interconnection})
and  (\ref{eq:gamma_0}) that
\begin{equation}\label{eq:T_i_defn}
\begin{split}
& T_i D_\gamma (\max_{\tau\in[t_{i-1},t_i]}\|\bfx(\tau)+\bfd(\tau)\|_{\varepsilon} + M h(t_0)\sigma_{i-1})\\
&\geq h(t_0)(\sigma_{i-1}-\sigma_{i}).
\end{split}
\end{equation}
Since
\[
\max_{\tau\in[t_{i-1},t_i]}\|\bfx(\tau)+\bfd(\tau)\|_{\varepsilon}=\|\bfx(\tau)+\bfd(\tau)\|_{\infty,[t_{i-1},t_i]}-\varepsilon
\]
if
$\|\bfx(\tau)+\bfd(\tau)\|_{\infty,[t_{i-1},t_i]}>\varepsilon$,
and
$\max_{\tau\in[t_{i-1},t_i]}\|\bfx(\tau)+\bfd(\tau)\|_{\varepsilon}=0$
overwise, we can see from (\ref{eq:T_i_defn}) that
\begin{equation}\label{eq:T_i_defn:1}
T_i\geq \left\{ \begin{array}{l} \frac{h(t_0)(\sigma_{i-1}-\sigma_i)}{D_\gamma(\|\bfx(\tau)+\bfd(\tau)\|_{\infty,[t_{i-1},t_i]}-\varepsilon+M h(t_0)\sigma_{i-1})},\\ \ \  \ \ \ \ | \|\bfx(\tau)+\bfd(\tau)\|_{\infty,[t_{i-1},t_i]}>\varepsilon;\\
\frac{\sigma_{i-1}-\sigma_i}{D_\gamma M \sigma_{i-1}},  \   \|\bfx(\tau)+\bfd(\tau)\|_{\infty,[t_{i-1},t_i]}\leq\varepsilon.
\end{array}\right.
\end{equation}
Pick
\begin{equation}\label{eq:tau_choice}
 \tau^\ast=\beta^{-1}\left({d}/{\kappa}\right), \ d\in(0,1),
\end{equation}
and select the value of $D_\gamma$ such that
(\ref{eq:gamma_choice}) holds. Consider two cases: a) $\|\bfx(\tau)+\bfd(\tau)\|_{\infty,[t_{i-1},t_i]}-\varepsilon\leq 0$ and b)
$\|\bfx(\tau)+\bfd(\tau)\|_{\infty,[t_{i-1},t_i]}-\varepsilon>0$. With respect to case a) we immediately observe that $T_i\geq \tau^\ast$ for all $i$.

Consider case b). Given that
$\|\bfx(\tau)+\bfd(\tau)\|_{\infty,[t_{0},t_1]}-\varepsilon\leq
\beta(0)\|\bfx(t_0)\|+c h(t_0)+\Delta+\Delta_d-\varepsilon$,
conditions (\ref{eq:gamma_choice}), (\ref{eq:epsilon_choice}),
and (\ref{eq:tau_choice}) imply $D_\gamma\leq
\frac{\kappa-1}{\kappa}\frac{h(t_0)}{\beta(0)\|\bfx(t_0)\|+(c+M)
|h(t_0)|}\frac{1}{\tau^\ast} \leq
\frac{h(t_0)(\sigma_0-\sigma_1)}{(\|\bfx(\tau)+\bfd(\tau)\|_{\infty,[t_{0},t_1]}-\varepsilon+M
|h(t_0)|)}\frac{1}{\tau^\ast}$. This, as follows from
(\ref{eq:T_i_defn:1}), guarantees that $T_1\geq \tau^\ast$.
Suppose that there is  an $i\geq 2$
such that $T_j\geq \tau^\ast$ for all $1\leq j\leq i-1$. We
will now show that the following implication  holds
$T_{i-1}\geq \tau^\ast\Rightarrow T_{i}\geq \tau^\ast$. This
will ensure that (\ref{eq:limit_h}) is satisfied and,
consequently, that the lemma hold. Consider
$\|\bfx(\tau)\|_{\infty,[t_{i-1},t_i]}$;
(\ref{eq:interconnection}) and (\ref{eq:h_i}) imply that:
$\|\bfx(\tau)\|_{\infty,[t_{i-1},t_i]}\leq
\beta(0)\|\bfx(t_{i-1})\|+c\sigma_{i-1}h(t_0)+\Delta$.
Estimating $\|\bfx(t_{i-1})\|$ from above, according to
(\ref{eq:interconnection}), results in
\[
\begin{split}
&\|\bfx(\tau)\|_{\infty,[t_{i-1},t_i]}\leq\beta(0)[\beta(T_{i-1})\|\bfx(t_{i-2})\|+c\sigma_{i-2}h(t_0)]\\
&+\beta(0)\Delta+ch(t_0)\sigma_{i-1}+\Delta\leq \beta(0)\beta(\tau^\ast) \|\bfx(t_{i-2})\| + P_1,
\end{split}
\]
where $P_1=\beta(0)c
\sigma_{i-2}h(t_0)+c\sigma_{i-1}h(t_0)+\beta(0)\Delta+\Delta$.
Invoking (\ref{eq:interconnection}) in order to express an
upper bound for $\|\bfx(t_{i-2})\|$ in terms of
$\|\bfx(t_{i-3})\|$  leads to
\[
\begin{split}
&\|\bfx(\tau)\|_{\infty,[t_{i-1},t_i]}\leq\beta(0)\beta^2(\tau^\ast)\|\bfx(t_{i-3})\|+P_2,
\end{split}
\]
where $P_2=c
h(t_0)\beta(0)[\beta(\tau^\ast)\sigma_{i-3}+\sigma_{i-2}]+c\sigma_{i-1}h(t_0)+
 \beta(0)[\beta(\tau^\ast)\Delta+\Delta]+\Delta$, and
\[
\begin{split}
&\|\bfx(\tau)\|_{\infty,[t_{i-1},t_i]}\leq\beta(0)\beta^3(\tau^\ast)\|\bfx(t_{i-4})\|+P_3,
\end{split}
\]
where
\[
\begin{array}{l}
P_3=c h(t_0)\beta(0)[\beta^2(\tau^\ast)\sigma_{i-4}+\beta(\tau^\ast)\sigma_{i-3}+\sigma_{i-2}]\\
+c\sigma_{i-1}h(t_0)+ \Delta\beta(0)[\beta(\tau^\ast)^2+\beta(\tau^\ast)+1]+\Delta\\
=c
h(t_0)\beta(0)\big[\sum_{j=0}^2\beta^j(\tau^\ast)\sigma_{i-j-2}\big]+ch(t_0)\sigma_{i-1}\\
+\Delta
\beta(0)\big[\sum_{j=0}^2\beta^j(\tau^\ast)\big]+\Delta.
\end{array}
\]
After $i-1$ steps we obtain
\begin{equation}\label{eq:x_i_est}
\begin{split}
\|\bfx(\tau)\|_{\infty,[t_{i-1},t_i]}&\leq \beta(0)\beta^{i-1}(\tau^\ast)\|\bfx(t_0)\|+P_{i-1},
\end{split}
\end{equation}
where
\[
\begin{array}{l}
P_{i-1}=ch(t_0)\beta(0)\big[\sum_{j=0}^{i-2}\beta^j(\tau^\ast)\sigma_{i-j-2}\big]+ ch(t_0)\sigma_{i-1}\\
+\Delta
\beta(0)\big[\sum_{j=0}^{i-2}\beta^j(\tau^\ast)\big]+\Delta.
\end{array}
\]
The values of $T_i$, as follows from (\ref{eq:T_i_defn:1}), satisfy:
\begin{eqnarray}\label{eq:T_i_bound}
T_i&\geq& \frac{\sigma_{i-1}-\sigma_{i}}{\sigma_{i-1}}h(t_0)\big[D_{\gamma}\sigma_{i-1}^{-1}\big(\|\bfx(\tau)+\bfd(\tau)\|_{\infty,[t_{i-1},t_i]}-\varepsilon\nonumber\\
& & +Mh(t_0)\sigma_{i-1}\big)\big]^{-1}.
\end{eqnarray}
Consider
\[
\sigma_{i-1}^{-1}\left(\|\bfx(\tau)+\bfd(\tau)\|_{\infty,[t_{i-1},t_i]}-\varepsilon+M
h(t_0)\sigma_{i-1}\right).
\]
Taking (\ref{eq:x_i_est}) into
account we derive that:
\[
\begin{array}{l}
\sigma_{i-1}^{-1}\left(\|\bfx(\tau)+\bfd(\tau)\|_{\infty,[t_{i-1},t_i]}-\varepsilon+M h(t_0)\sigma_{i-1}\right)\\
\leq \beta(0) \beta^{i-1}(\tau^\ast)\kappa^{i-1}\|\bfx(t_0)\|+k^{i-1}P_{i-1}+k^{i-1}\Delta_d\\
-k^{i-1}\varepsilon + M h(t_0)= \beta(0) \beta^{i-1}(\tau^\ast)\kappa^{i-1}\|\bfx(t_0)\|\\
+ch(t_0)\beta(0)\kappa\big[\sum_{j=0}^{i-2}\beta^j(\tau^\ast)\kappa^j\big]+(c+M)h(t_0)\\
+\kappa^{i-1}\big[\Delta\big(\beta(0)\sum_{j=0}^{i-2}\beta^j(\tau^\ast)+1\big)+\Delta_d-\varepsilon\big].
\end{array}
\]
Noticing that $\tau^\ast$ is chosen in accordance with
(\ref{eq:tau_choice}) one can therefore obtain:
\[
\begin{array}{l}
\sigma_{i-1}^{-1}\left(\|\bfx(\tau)+\bfd(\tau)\|_{\infty,[t_{i-1},t_i]}-\varepsilon+Mh(t_0)\sigma_{i-1}\right)\\
\leq \beta(0)\|\bfx(t_0)\|+(c+M)h(t_0)+ch(t_0)\beta(0)\kappa\sum_{j=0}^{i-2}d^j+\\
\kappa^{i-1}\big[\Delta\big(\beta(0)\sum_{j=0}^{i-2}\frac{d}{\kappa}^j+1\big)+\Delta_d-\varepsilon\big]\\
\leq \beta(0)\|\bfx(t_0)\|+h(t_0)\big(c\big(1+\frac{\beta(0)\kappa}{1-d} \big)+M\big)+\\
k^{i-1}\big[\Delta\big(\frac{\beta(0)}{1-d/k}+1\big)+\Delta_d-\varepsilon\big].
\end{array}
\]
Condition (\ref{eq:epsilon_choice}) implies that
$\Delta\left(\frac{\beta(0)}{1-d/k}+1\right)+\Delta_d-\varepsilon\leq
0$. Hence
$\sigma_{i-1}^{-1}\left(\|\bfx(\tau)+\bfd(\tau)\|_{\infty,[t_{i-1},t_i]}-\varepsilon+Mh(t_0)\sigma_{i-1}\right)\leq
\beta(0)\|\bfx(t_0)\|+h(t_0)\left(c\left(1+\frac{\beta(0)\kappa}{1-d}
\right)+M\right)$. Substituting the latter estimate into
(\ref{eq:T_i_bound}) and using (\ref{eq:gamma_choice}) yields
$T_{i}\geq  \tau^\ast$ as required. Therefore $T_i\geq
\tau^\ast$ for all $i\geq 1$, and (\ref{eq:limit_h}) holds.
Thus the trajectory $h(t)$ is bounded from above and below, and
hence so is the trajectory $\bfx(t)$. $\square$

\end{document}